\documentclass[11pt]{amsart}
\usepackage{amscd,amssymb,amsthm}
\usepackage{graphicx}

\newtheorem{theorem}{Theorem}[section]
\newtheorem{remark}{Remark}[section]

\newtheorem{proposition}{Proposition}[section]
\numberwithin{equation}{section}

\title[To appear in: \textit{Journal of Concrete and Applicable
Mathematics}]{Bounds On Triangular Discrimination, Harmonic Mean
and Symmetric Chi-square Divergences}

\author{Inder Jeet Taneja}
\address{Departamento de Matem\'{a}tica\\ Universidade Federal
de Santa Catarina\\
88.040-900 Florian\'{o}polis, SC, Brazil}

\email{taneja@mtm.ufsc.br}

\urladdr{http://www.mtm.ufsc.br/$\sim $taneja}


\thanks{To appear in: \textit{Journal of Concrete and Applicable
Mathematics (2005)}}


\begin{document}

\begin{abstract}
There are many information and divergence measures exist in the
literature on information theory and statistics. The most famous
among them are Kullback-Leiber \cite{kul} \textit{relative
information} and Jeffreys \cite{jef} \textit{J-divergence}. The
measures like \textit{Bhattacharya distance}, \textit{Hellinger
discrimination}, \textit{Chi-square divergence},
\textit{triangular discrimination} and \textit{harmonic mean
divergence} are also famous in the literature on statistics. In
this paper we have obtained bounds on \textit{triangular
discrimination} and \textit{symmetric chi-square divergence} in
terms of \textit{relative information of type s} using
Csisz\'{a}r's \textit{f-divergence}. A relationship among
\textit{triangular discrimination} and \textit{harmonic mean
divergence} is also given.
\end{abstract}

\maketitle

\section{Introduction}

Let
\[
\Gamma _n = \left\{ {P = (p_1 ,p_2 ,...,p_n )\left| {p_i > 0,\sum\limits_{i
= 1}^n {p_i = 1} } \right.} \right\},
\quad
n \geqslant 2,
\]

\noindent be the set of all complete finite discrete probability
distributions. For all $P,Q \in \Gamma _n $, the following
measures are well known in the literature on information theory
and statistics:\\

\noindent \textbf{$\bullet$ Bhattacharya Distance} (Bhattacharya
\cite{bha})
\begin{equation}
\label{eq1} B(P\vert \vert Q) = \sum\limits_{i = 1}^n \sqrt {p_i
q_i }.
\end{equation}

\noindent \textbf{$\bullet$ Hellinger discrimination} (Hellinger
\cite{hel})
\begin{equation}
\label{eq2}
h(P\vert \vert Q) = 1 - B(P\vert \vert Q) = \frac{1}{2}\sum\limits_{i = 1}^n
{(\sqrt {p_i } - \sqrt {q_i } )^2} .
\end{equation}

\noindent \textbf{$\bullet$ $\chi ^2 - $Divergence} (Pearson
\cite{pea})
\begin{equation}
\label{eq3}
\chi ^2(P\vert \vert Q) = \sum\limits_{i = 1}^n {\frac{(p_i - q_i )^2}{q_i
}} = \sum\limits_{i = 1}^n {\frac{p_i^2 }{q_i } - 1} .
\end{equation}

\noindent \textbf{$\bullet$ Relative Information} (Kullback and
Leibler \cite{kul})
\begin{equation}
\label{eq4} K(P\vert \vert Q) = \sum\limits_{i = 1}^n {p_i \ln
(\frac{p_i }{q_i })}.
\end{equation}

The above four measures can be obtained as particular or limiting
case of the \textit{relative information of type s}. This
measure is given by\\

\noindent \textbf{$\bullet$ Relative Information of Type s}
\begin{equation}
\label{eq5} \Phi _s (P\vert \vert Q) = \begin{cases}
 {K_s (P\vert \vert Q) = \left[ {s(s - 1)} \right]^{ - 1}\left[
{\sum\limits_{i = 1}^n {p_i^s q_i^{1 - s} - 1} } \right],} & {s \ne 0,1} \\
 {K(Q\vert \vert P) = \sum\limits_{i = 1}^n {q_i \ln \left( {\frac{q_i }{p_i
}} \right)} ,} & {s = 0} \\
 {K(P\vert \vert Q) = \sum\limits_{i = 1}^n {p_i \ln \left( {\frac{p_i }{q_i
}} \right)} ,} & {s = 1} \\
\end{cases}.
\end{equation}

The measure (\ref{eq5}) admits the following interesting
particular cases:
\begin{itemize}
\item[(i)] $ \Phi _{ - 1} (P\vert \vert Q) = \frac{1}{2}\chi
^2(Q\vert \vert P)$.\\
\item[(ii)] $ \Phi _0 (P\vert \vert Q) = K(Q\vert \vert P)$.\\
\item[(iii)] $ \Phi _{1 / 2} (P\vert \vert Q) = 4\left[ {1 -
B(P\vert \vert Q)} \right] = 4h(P\vert \vert Q)$.\\
\item[(iv)] $ \Phi _1 (P\vert \vert Q) = K(P\vert \vert Q)$.\\
\item[(v)] $ \Phi _2 (P\vert \vert Q) = \frac{1}{2}\chi ^2(P\vert
\vert Q)$.\\
\end{itemize}

Thus we observe that $\Phi _2 (P\vert \vert Q) = \Phi _{ - 1}
(Q\vert \vert P)$ and $\Phi _1 (P\vert \vert Q) = \Phi _0 (Q\vert
\vert P)$.\\

For more studies on the measure (\ref{eq5}) refer to Liese and
Vajda \cite{liv}, Vajda \cite{vaj}, Taneja \cite{tan1},
\cite{tan2}, \cite{tan4} and Cerone et al. \cite{cdo}.\\

Recently Taneja \cite{tan5} and Taneja and Kumar \cite{tak}
studied the (\ref{eq5}) and obtained bounds in terms of the
measures (\ref{eq1})-(\ref{eq4}). Here we shall extend the our
study for the other measures known in the literature as
\textit{triangular discrimination}, \textit{harmonic mean
divergence} and \textit{symmetric chi-square
divergence}.\\

The \textit{triangular discrimination} is given by
\begin{equation}
\label{eq6}
\Delta (P\vert \vert Q) = \sum\limits_{i = 1}^n {\frac{(p_i - q_i )^2}{p_i +
q_i }} .
\end{equation}

After simplification, we can write
\begin{equation}
\label{eq7}
\Delta (P\vert \vert Q) = 2\left[ {1 - W(P\vert \vert Q)} \right],
\end{equation}

\noindent where
\begin{equation}
\label{eq8}
W(P\vert \vert Q) = \sum\limits_{i = 1}^n {\frac{2p_i q_i }{p_i + q_i }} ,
\end{equation}

\noindent is the well known \textit{harmonic mean divergence}.\\

We observe that the measures (\ref{eq3}) and (\ref{eq4}) are not
symmetric with respect to probability distributions. The symmetric
version of the measure (\ref{eq4}) famous as
Jeffreys-Kullback-Leiber \textit{J-divergence} is given by
\begin{equation}
\label{eq9}
J(P\vert \vert Q) = K(P\vert \vert Q) + K(Q\vert \vert P).
\end{equation}

Let us consider the \textit{symmetric chi-square divergence} given
by
\begin{equation}
\label{eq10} \Psi(P\vert \vert Q) = \chi ^2(P\vert \vert Q) + \chi
^2(Q\vert \vert P)
 = \sum\limits_{i = 1}^n {\frac{(p_i - q_i )^2(p_i + q_i )}{p_i q_i }} .
\end{equation}

Dragomir \cite{dsb} studied the measure (\ref{eq10}) and obtained
interesting result relating it to \textit{triangular
discrimination} and \textit{J-divergence}.

Some studies on the measures (\ref{eq6}) and (\ref{eq8}) can be
seen in Dragomir \cite{dra5}, \cite{dra6} and Tops$\o$e
\cite{top}. Recently, Taneja \cite{tan5} and Taneja and Kumar
\cite{tak} studied the measure (\ref{eq5}) and obtained bounds in
terms of the measures (\ref{eq1})-(\ref{eq4}). Similar kind of
bounds on the measure (\ref{eq9}) are recently obtained by Taneja
\cite{tan6}. In this paper, we shall extend the our study for
\textit{triangular discrimination} and \textit{symmetric
chi-square divergence}. In order to obtain bounds on these
measures we make use of Csisz\'{a}r's \cite{csi}
\textit{f-divergence}.

\section{Csisz\'{a}r's $f-$Divergence and Its Particular Cases}

Given a convex function $f:[0,\infty ) \to \mathbb{R}$, the $f -
$divergence measure introduced by Csisz\'{a}r \cite{csi} is given
by
\begin{equation}
\label{eq12}
C_f (P\vert \vert Q) =
\sum\limits_{i = 1}^n {q_i f\left( {\frac{p_i }{q_i }} \right)} ,
\end{equation}

\noindent where $P,Q \in \Gamma _n $.\\

It is well known in the literature \cite{csi} that \textit{if $f$
is convex and normalized, i.e., $f(1) = 0$, then the Csisz\'{a}r's
function $C_f (P\vert \vert Q)$ is nonnegative and convex in the
pair of probability distribution $(P,Q) \in \Gamma _n \times
\Gamma _n $}.

Here below we shall give the measures (\ref{eq6}) and
(\ref{eq10}) being examples of the measure (\ref{eq12}).\\

\textbf{Example 2.1. }(\textit{Triangular discrimination}). Let us
consider
\begin{equation}
\label{eq13}
f_\Delta (x) = \frac{(x - 1)^2}{x + 1},
\quad
x \in (0,\infty )
\end{equation}

\noindent in (\ref{eq12}), we have
\[
C_f (P\vert \vert Q) = \Delta (P\vert \vert Q) = \sum\limits_{i = 1}^n
{\frac{(p_i - q_i )^2}{p_i + q_i }} ,
\]

\noindent where $\Delta (P\vert \vert Q)$ is as given by
(\ref{eq6}).

Moreover,
\begin{equation}
\label{eq14}
{f}'_\Delta (x) = \frac{(x - 1)(x + 3)}{(x + 1)^2}
\end{equation}

\noindent and
\begin{equation}
\label{eq15}
{f}''_\Delta (x) = \frac{8}{(x + 1)^3}.
\end{equation}

Thus we have ${f}''_\Delta (x) > 0$ for all $x > 0$, and hence,
$f_\Delta (x)$ is \textit{strictly convex} for all $x > 0$. Also,
we have $f_\Delta (1) = 0$. In view of this we can say that the
\textit{triangular discrimination} is \textit{nonnegative} and
\textit{convex} in the pair of probability distributions $(P,Q)
\in \Gamma _n \times \Gamma _n $.\\

\textbf{Example 2.2. }(\textit{Symmetric chi-square divergence}).
Let us consider
\begin{equation}
\label{eq16} f_\Psi (x) = \frac{(x - 1)^2(x + 1)}{x}, \quad x \in
(0,\infty )
\end{equation}

\noindent in (\ref{eq12}), we have
\[
C_f (P\vert \vert Q) = \Psi(P\vert \vert Q) = \sum\limits_{i =
1}^n {\frac{(p_i - q_i )^2(p_i + q_i )}{p_i q_i }} ,
\]

\noindent where $\Psi\left( {P\vert \vert Q} \right)$ is as given
by (\ref{eq10}).

Moreover,
\begin{equation}
\label{eq17} {f}'_\Psi (x) = \frac{(x - 1)(2x^2 + x + 1)}{x^2}
\end{equation}

\noindent and
\begin{equation}
\label{eq18} {f}''_\Psi (x) = \frac{2(x^3 + 1)}{x^3}.
\end{equation}

Thus we have ${f}''_\Psi (x) > 0$ for all $x > 0$, and hence,
$f_\Psi (x)$ is \textit{strictly convex} for all $x > 0$. Also, we
have $f_\Psi (1) = 0$. In view of this we can say that the
\textit{symmetric chi-square divergence} is \textit{nonnegative}
and \textit{convex} in the pair of probability distributions
$(P,Q) \in \Gamma _n \times \Gamma _n $.

\section{Csisz\'{a}r's $f-$Divergence and Relative
Information of Type $s$}

During past years Dragomir done a lot of work giving bounds on
Csisz\'{a}r's \textit{$f - $divergence}. Here below we shall
summarize the some his results \cite{dra1}, \cite{dra2},
\cite{dra4}.

\begin{theorem} \label{the31} Let $f:\mathbb{R}_ + \to [0,\infty )$ be
differentiable convex and normalized i.e., $f(1) = 0$. If $P,Q \in
\Gamma _n $, then we have
\begin{equation}
\label{eq19}
0 \leqslant C_f (P\vert \vert Q) \leqslant \rho _{C_f } (P\vert \vert Q),
\end{equation}

\noindent where $\rho _{C_f } (P\vert \vert Q)$ is given by
\begin{equation}
\label{eq20} \rho _{C_f } (P\vert \vert Q) = C_{f'} \left(
{\frac{P^2}{Q}\vert \vert P} \right) - C_{f'} (P\vert \vert Q) =
\sum\limits_{i = 1}^n {(p_i - q_i )} {f}'(\frac{p_i }{q_i }).
\end{equation}

If $P,Q \in \Gamma _n $ are such that
\[
0 < r \leqslant \frac{p_i }{q_i } \leqslant R < \infty ,
\quad
\forall i \in \{1,2,...,n\},
\]

\noindent for some $r$ and $R$ with $0 < r \leqslant 1 \leqslant R
< \infty $, then we have the following inequalities:
\begin{equation}
\label{eq21}
0 \leqslant C_f (P\vert \vert Q) \leqslant \alpha _{C_f } (r,R),
\end{equation}
\begin{equation}
\label{eq22}
0 \leqslant C_f (P\vert \vert Q) \leqslant \beta _{C_f } (r,R)
\end{equation}

\noindent and
\begin{align}
0 & \leqslant \beta _{C_f } (r,R) - C_f (P\vert \vert
Q)\label{eq23}\\
& \leqslant \gamma _{C_f } (r,R)\left[ {(R - 1)(1 - r) - \chi
^2(P\vert \vert Q)} \right] \leqslant \alpha _{C_f } (r,R),\notag
\end{align}

\noindent where
\begin{align}
\alpha _{C_f } (r,R) & = \frac{1}{4}(R - r)^2\mbox{ }\gamma _{C_f
}
(r,R),\label{eq24}\\
 \beta _{C_f } (r,R) & = \frac{(R - 1)f(r) + (1 -
r)f(R)}{R - r}\label{eq25}\\
\intertext{and}
\gamma _{C_f } (r,R) &= \frac{{f}'(R) - {f}'(r)}{R
- r}.\label{eq26}
\end{align}
\end{theorem}

The following proposition is due to Taneja \cite{tan5} and Taneja
and Kumar \cite{tak} and is a consequence of the above theorem.

\begin{proposition} \label{pro31} Let $P,Q \in \Gamma _n $ and $s \in
\mathbb{R}$, then we have
\begin{equation}
\label{eq27}
0 \leqslant
\Phi _s (P\vert \vert Q)
 \leqslant \rho _{\Phi _s} (P\vert \vert Q),
\end{equation}

\noindent where
\begin{align}
\rho _{\Phi _s} (P\vert \vert Q) & = C_{\phi _s ^\prime } \left(
{\frac{P^2}{Q}\vert \vert P} \right) - C_{\phi _s ^\prime } \left(
{P\vert \vert Q} \right)\label{eq28}\\
& = \begin{cases}
 {(s - 1)^{ - 1}\sum\limits_{i = 1}^n {(p_i - q_i )\left( {\frac{p_i }{q_i
}} \right)^{s - 1},} } & {s \ne 1} \\
 {\sum\limits_{i = 1}^n {(p_i - q_i )\ln \left( {\frac{p_i }{q_i }}
\right),} } & {s = 1} \\
\end{cases}.\notag
\end{align}

If there exists $r,R$ $(0 < r \leqslant 1 \leqslant R < \infty )$
such that
\[
0 < r \leqslant \frac{p_i }{q_i } \leqslant R < \infty ,
\quad
\forall i \in \{1,2,...,n\},
\]

\noindent then we have the following inequalities
\begin{equation}
\label{eq29} 0 \leqslant \Phi _s (P\vert \vert Q) \leqslant \alpha
_{\Phi _s} (r,R),
\end{equation}
\begin{equation}
\label{eq30} 0 \leqslant \Phi _s (P\vert \vert Q) \leqslant \beta
_{\Phi _s} (r,R)
\end{equation}

\noindent and
\begin{align}
0 & \leqslant \beta _{\Phi _s} (r,R) - \Phi _s (P\vert \vert
Q)\label{eq31}\\
& \leqslant \gamma _{\Phi _s} (r,R)\left[ {(R - 1)(1 - r) - \chi
^2(P\vert \vert Q)} \right] \leqslant \alpha _{\Phi _s}
(r,R),\notag
\end{align}

\noindent where
\begin{align}
\alpha _{\Phi _s} (r,R) & = \frac{1}{4}(R - r)^2\mbox{ }\gamma
_{\Phi _s } (r,R),\label{eq32}\\
\beta _{\Phi _s} (r,R) & = \begin{cases}
 {\frac{(R - 1)(r^s - 1) + (1 - r)(R^s - 1)}{(R - r)s(s - 1)},} & {s \ne
0,1} \\
 {\frac{(R - 1)\ln \frac{1}{r} + (1 - r)\ln \frac{1}{R}}{(R - r)},} & {s =
0} \\
 {\frac{(R - 1)r\ln r + (1 - r)R\ln R}{(R - r)},} & {s = 1} \\
\end{cases}\label{eq33}\\
\intertext{and} \gamma _{\Phi _s } (r,R) & =
\begin{cases}
 {\frac{R^{s - 1} - r^{s - 1}}{(R - r)(s - 1)},} & {s \ne 1} \\
 {\frac{\ln R - \ln r}{R - r},} & {s = 1} \\
\end{cases},\label{eq34}
\end{align}
\end{proposition}

We can also write $\gamma _{\Phi _s} (r,R)$ as follows
\begin{equation}
\label{eq35} \gamma _{\Phi _s} (r,R) = \begin{cases}
 {L_{s - 2}^{s - 2} (r,R),} & {s \ne 1} \\
 {L_{ - 1}^{ - 1} (r,R)} & {s = 1} \\
\end{cases},
\end{equation}

\noindent where $L_p (a,b)$ is the famous (Bullen, Mitrinovi\'{c}
and Vasi\'{c} \cite{bmv}) \textit{p-logarithmic power mean} given
by
\begin{equation}
\label{eq36} L_p (a,b) = \begin{cases}
 {\left[ {\frac{b^{p + 1} - a^{p + 1}}{(p + 1)(b - a)}}
\right]^{\frac{1}{p}},} & {p \ne - 1,0} \\
 {\frac{b - a}{\ln b - \ln a},} & {p = - 1} \\
 {\frac{1}{e}\left[ {\frac{b^b}{a^a}} \right]^{\frac{1}{b - a}},} & {p = 0}
\\
\end{cases},
\end{equation}

\noindent for all $p \in \mathbb{R}$, $a \ne b$.

The expression (\ref{eq28}) admits the following particular cases:
\begin{itemize}
\item[(i)] $\rho _{\Phi _{ - 1} } (P\vert \vert Q) = 3\,\Phi _3
(Q\vert \vert P) - \frac{1}{2}\chi ^2(Q\vert \vert P)$.\\

\item[(ii)] $\rho _{\Phi _0 } (P\vert \vert Q) = \chi ^2(Q\vert
\vert P)$.\\

\item[(iii)] $\rho _{\Phi _1 } (P\vert \vert Q) = J(P\vert \vert
Q)$.\\

\item[(iv)] $\rho _{\Phi _{1 / 2} } (P\vert \vert Q) =
2\sum\limits_{i = 1}^n {(q_i - p_i )\sqrt {\frac{q_i }{p_i }} }
$\\

\item[(v)] $\rho _{\Phi _2 } (P\vert \vert Q) = \chi ^2(P\vert
\vert Q)$\\
\end{itemize}

The expression (\ref{eq33}) admits the following particular cases:

\begin{itemize}
\item[(i)] $\beta _{\Phi _{ - 1} } (P\vert \vert Q) = \frac{(R -
1)(1 - r)}{2rR}$.\\

\item[(ii)] $\beta _{\Phi _0 } (r,R) = \frac{(R - 1)\ln
\frac{1}{r} + (1 - r)\ln \frac{1}{R}}{R - r}$.\\

\item[(iii)] $\beta _{\Phi _1 } (r,R) = \frac{(R - 1)r\ln r + (1 -
r)R\ln R}{R - r}$.\\

\item[(iv)] $\beta _{\Phi _{1 / 2} } (P\vert \vert Q) =
\frac{4(\sqrt R - 1)(1 - \sqrt r )}{\sqrt R + \sqrt r }$.\\

\item[(v)] $\beta _{\Phi _2 } (P\vert \vert Q) = \frac{(R - 1)(1 -
r)}{2}$.\\
\end{itemize}

The following theorem is due to Taneja \cite{tan5} and Taneja and
Kumar \cite{tak}.

\begin{theorem} \label{the32} Let $f:I \subset \mathbb{R}_ + \to
[0,\infty )$ the generating mapping is normalized, i.e., $f(1) =
0$ and satisfy the assumptions:

(i) $f$ is twice differentiable on $(r,R)$, where $0 \leqslant r
\leqslant 1 \leqslant R \leqslant \infty $;

(ii) there exists real constants $m,M$ such that $m < M$ and
\begin{equation}
\label{eq37}
m \leqslant x^{2 - s}{f}''(x) \leqslant M,
\quad
\forall x \in (r,R),
\quad
s \in \mathbb{R}.
\end{equation}

If $P,Q \in \Gamma _n $ are discrete probability distributions satisfying
the assumption
\[
0 < r \leqslant \frac{p_i }{q_i } \leqslant R < \infty ,
\]

\noindent then we have the inequalities:
\begin{equation}
\label{eq38}
m
\Phi _s (P\vert \vert Q)
 \leqslant C_f (P\vert \vert Q) \leqslant M
\Phi _s (P\vert \vert Q),
\end{equation}
\begin{align}
m & \left[ {\rho _{\Phi _s } (P\vert \vert Q) - \Phi _s (P\vert
\vert Q)} \right] \label{eq39} \\
& \leqslant \rho _{C_f } (P\vert \vert Q) - C_f (P\vert \vert
Q)\notag\\
& \leqslant M\left[ {\rho _{\Phi _s } (P\vert \vert Q) - \Phi _s
(P\vert \vert Q)} \right]\notag
\end{align}

\noindent and
\begin{align}
m & \left[ {\beta _{\Phi _s } (r,R) - \Phi _s (P\vert \vert Q)}
\right] \label{eq40}\\
& \leqslant \beta _{C_f } (r,R) - C_f (P\vert \vert Q)\notag\\
& \leqslant M\left[ {\beta _{\Phi _s } (r,R) - \Phi _s (P\vert
\vert Q)} \right],\notag
\end{align}

\noindent where $C_f (P\vert \vert Q)$, $\Phi _s (P\vert \vert
Q)$, $\rho _{C_f } (P\vert \vert Q), \rho _{\Phi _s } (P\vert
\vert Q)$, $\beta _{C_f } (r,R)$ and $\beta _{\Phi _s } (r,R)$ are
as given by (\ref{eq12}), (\ref{eq5}), (\ref{eq20}), (\ref{eq28}),
(\ref{eq25}) and (\ref{eq33}) respectively.
\end{theorem}

The above theorem unifies some of the results studied by Dragomir
\cite{dra3}, \cite{dra5}, \cite{dra6}.\\

In the papers Taneja \cite{tan5} and Taneja and Kumar \cite{tak}
considered the particular values of $s$ and $\Phi _s $ by taking
$s = -1$, $s = 0$, $s = \frac{1}{2}$, $s = 1$ and $s = 2$. The aim
here is to obtain results by taking different values of $f$ given
by examples 2.1-2.2, and then obtain particular cases for
different values of $s$.

\begin{remark} If is it not specified, from now onwards, it
is understood that, if there are $r,R$ then $0 < r \leqslant
\frac{p_i }{q_i } \leqslant R < \infty $, $\forall i \in
\{1,2,...,n\}$, with $0 < r \leqslant 1 \leqslant R < \infty $
where $P = (p_1 ,p_2 ,....,p_n ) \in \Gamma _n $ and $P = (q_1
,q_2 ,....,q_n ) \in \Gamma _n $.
\end{remark}

\section{Triangular Discrimination and Inequalities}

In this section, we shall give bounds on \textit{triangular
discrimination} based on the Theorems \ref{the31} and \ref{the32}.
\begin{theorem} \label{the41} For all $P,Q \in \Gamma _n $, we have
the following inequalities
\begin{equation}
\label{eq41}
0 \leqslant \Delta (P\vert \vert Q) \leqslant \rho _\Delta (P\vert \vert
Q),
\end{equation}

\noindent where
\begin{equation}
\label{eq42}
\rho _\Delta (P\vert \vert Q) = \sum\limits_{i = 1}^n {\left( {\frac{p_i -
q_i }{p_i + q_i }} \right)^2(p_i + 3q_i )} .
\end{equation}

If there exists $r,R$ $(0 < r \leqslant 1 \leqslant R < \infty )$
such that
\[
0 < r \leqslant \frac{p_i }{q_i } \leqslant R < \infty ,
\quad
\forall i \in \{1,2,...,n\},
\]

\noindent then we have the following inequalities:
\begin{equation}
\label{eq43}
0 \leqslant \Delta (P\vert \vert Q) \leqslant \alpha _\Delta (r,R),
\end{equation}
\begin{equation}
\label{eq44}
0 \leqslant \Delta (P\vert \vert Q) \leqslant \beta _\Delta (r,R)
\end{equation}

\noindent and
\begin{align}
0 & \leqslant \beta _\Delta (r,R) - \Delta (P\vert \vert
Q)\label{eq45}\\
& \leqslant \gamma _\Delta (r,R)\left[ {(R - 1)(1 - r) - \chi
^2(P\vert \vert Q)} \right]
 \leqslant \alpha _\Delta (r,R),\notag
\end{align}

\noindent where
\begin{align}
\alpha _\Delta (r,R) & = \frac{1}{4}(R - r)\left[ {\frac{(R - 1)(R
+ 3)}{(R + 1)^2} + \frac{(1 - r)(r + 3)}{(r + 1)^2}}
\right],\label{eq46}\\
\beta _\Delta (r,R) & = \frac{2(R - 1)(1 - r)}{(R + 1)(1 +
r)}\label{eq47} \\
\intertext{and} \gamma _\Delta (r,R)& = (R - r)^{ - 1}\left[
{\frac{(R - 1)(R + 3)}{(R + 1)^2} + \frac{(1 - r)(r + 3)}{(r +
1)^2}} \right].\label{eq48}
\end{align}
\end{theorem}

\begin{proof} Follows from the Theorem \ref{the31} by
considering $f$ by $f_\Delta $ and making necessary calculations.
\end{proof}

\begin{theorem} \label{the42} Let $P,Q \in \Gamma _n $ and $s \in
\mathbb{R}$. Let there exists $r,R$ $(0 < r \leqslant 1 \leqslant
R < \infty )$ such that $0 < r \leqslant \frac{p_i }{q_i }
\leqslant R < \infty $, $\forall i \in \{1,2,...,n\}$.

$(a)$ For $s \leqslant - 1$, we have the following inequalities:
\begin{equation}
\label{eq49}
\frac{8r^{2 - s}}{(r + 1)^3}\Phi _s (P\vert \vert Q)
 \leqslant \Delta (P\vert \vert Q) \leqslant
\frac{8R^{2 - s}}{(R + 1)^3}\Phi _s (P\vert \vert Q),
\end{equation}
\begin{align}
\frac{8r^{2 - s}}{(r + 1)^3}& \left[ {\rho_{\Phi_s} (P\vert \vert
Q) - \Phi _s (P\vert \vert Q)} \right] \label{eq50}\\
& \leqslant \Delta ^* (P\vert \vert Q)
 \leqslant \frac{8R^{2 - s}}{(R + 1)^3}\left[ {\rho_{\Phi_s} (P\vert \vert Q) -
\Phi _s (P\vert \vert Q)} \right]\notag
\end{align}

\noindent and
\begin{align}
\frac{8r^{2 - s}}{(r + 1)^3}& \left[ {\beta_{\Phi_s} (r,R) - \Phi
_s
(P\vert \vert Q)} \right] \label{eq51}\\
& \leqslant \beta _\Delta (r,R) - \Delta (P\vert \vert Q)\notag\\
& \leqslant \frac{8R^{2 - s}}{(R + 1)^3}\left[ {\beta_{\Phi_s}
(r,R) - \Phi _s (P\vert \vert Q)} \right].\notag
\end{align}

$(b)$ For $s \geqslant 2$, we have the following inequalities:
\begin{equation}
\label{eq52}
\frac{8R^{2 - s}}{(R + 1)^3}\Phi _s (P\vert \vert Q)
 \leqslant \Delta (P\vert \vert Q) \leqslant
\frac{8r^{2 - s}}{(r + 1)^3}\Phi _s (P\vert \vert Q),
\end{equation}
\begin{align}
\frac{8R^{2 - s}}{(R + 1)^3}& \left[ {\rho_{\Phi_s} (P\vert \vert
Q) -
\Phi _s (P\vert \vert Q)} \right] \label{eq53}\\
& \leqslant \Delta ^* (P\vert \vert Q)
 \leqslant \frac{8r^{2 - s}}{(r + 1)^3}\left[ {\rho_{\Phi_s} (P\vert \vert Q) -
\Phi _s (P\vert \vert Q)} \right]\notag
\end{align}

\noindent and
\begin{align}
\frac{R^{1 - s}}{(R + 1)^2}& \left[ {\beta_{\Phi_s} (r,R) - \Phi
_s (P\vert
\vert Q)} \right]\label{eq54}\\
&  \leqslant \beta _\Delta (r,R) - \Delta (P\vert \vert Q)\notag\\
& \leqslant \frac{r^{1 - s}}{(r + 1)^2}\left[ {\beta_{\Phi_s}
(r,R) - \Phi _s (P\vert \vert Q)} \right],\notag
\end{align}

\noindent where
\begin{equation}
\label{eq55}
\Delta ^ * (P\vert \vert Q) = \rho _\Delta (P\vert \vert Q) - \Delta (P\vert
\vert Q) = 2\sum\limits_{i = 1}^n {q_i \left( {\frac{p_i - q_i }{p_i + q_i
}} \right)} ^2.
\end{equation}
\end{theorem}

\begin{proof} Let us consider
\begin{equation}
\label{eq56}
g_\Delta (x) = x^{2 - s}{f}''_\Delta (x) = \frac{8x^{2 - s}}{(x + 1)^3},
\quad
x \in (0,\infty ),
\end{equation}

\noindent where ${f}''_\Delta (x)$ is as given by (\ref{eq15}).

We have
\begin{equation}
\label{eq57}
{g}'_\Delta (x) = - \frac{8x^{1 - s}\left[ {(s + 1)x + (s - 2)} \right]}{(x
+ 1)^4}
\begin{cases}
 { \geqslant 0,} & {s \leqslant - 1} \\
 { \leqslant 0,} & {s \geqslant 2} \\
\end{cases}.
\end{equation}

In view of (\ref{eq57}), we conclude the followings:
\begin{equation}
\label{eq58} m = \mathop {\inf }\limits_{x \in [r,R]} g(x) =
\mathop {\min }\limits_{x \in [r,R]} g(x) = \begin{cases}
 {\frac{8r^{2 - s}}{(r + 1)^3},} & {s \leqslant - 1} \\
 {\frac{8R^{2 - s}}{(R + 1)^3},} & {s \geqslant 2} \\
\end{cases}
\end{equation}

\noindent and
\begin{equation}
\label{eq59} M = \mathop {\sup }\limits_{x \in [r,R]} g(x) =
\mathop {\max }\limits_{x \in [r,R]} g(x) = \begin{cases}
 {\frac{8R^{2 - s}}{(R + 1)^3},} & {s \leqslant - 1} \\
 {\frac{8r^{2 - s}}{(r + 1)^3},} & {s \geqslant 2} \\
\end{cases}.
\end{equation}

From (\ref{eq58}) and (\ref{eq59}) and Theorem \ref{the32}, we
have the required proof.
\end{proof}

The following propositions are the particular cases of the above
theorem.

\begin{proposition} \label{pro41} We have the following bounds in
terms of $\chi ^2 - $divergence:
\begin{equation}
\label{eq60}
\frac{4r^3}{(r + 1)^3}\chi ^2(Q\vert \vert P) \leqslant
\Delta (P\vert \vert
Q) \leqslant \frac{4R^3}{(R + 1)^3}\chi ^2(Q\vert \vert P),
\end{equation}
\begin{align}
\frac{8 r^3}{(r + 1)^3} &  \left[3 \,\Phi _3 (Q\vert \vert P) -
\chi^2 (Q\vert \vert P) \right] \label{eq61}\\
& \leqslant \Delta ^ * (P\vert \vert Q)
 \leqslant \frac{8 R^3}{(R + 1)^3}
 \left[ {3\,\Phi _3 (Q\vert \vert P) - \chi
^2(Q\vert \vert P)} \right]\notag
\end{align}

\noindent and
\begin{align}
\frac{4r^3}{(r + 1)^3}& \left[ {\frac{(R - 1)(1 - r)}{rR} - \chi
^2(Q\vert \vert P)} \right] \label{eq62}\\
& \leqslant \frac{2(R - 1)(1 - r)}{(R + 1)(1 + r)} - \Delta
(P\vert \vert Q)\notag\\
& \leqslant \frac{4R^3}{(R + 1)^3}\left[ {\frac{(R - 1)(1 -
r)}{rR} - \chi ^2(Q\vert \vert P)} \right].\notag
\end{align}
\end{proposition}

\begin{proof} Take $s = - 1$ in (\ref{eq49}), (\ref{eq50})
and (\ref{eq51}) we get respectively (\ref{eq60}), (\ref{eq61})
and (\ref{eq62}).
\end{proof}

\begin{proposition} \label{pro42} We have the following bounds in
terms of $\chi ^2 - $divergence:
\begin{equation}
\label{eq63}
\frac{4}{(R + 1)^3}\chi ^2(P\vert \vert Q) \leqslant \Delta (P\vert \vert Q)
\leqslant \frac{4}{(r + 1)^3}\chi ^2(P\vert \vert Q),
\end{equation}
\begin{equation}
\label{eq64} \frac{4}{(R + 1)^3}\chi ^2(P\vert \vert Q) \leqslant
\Delta ^ * (P\vert \vert Q) \leqslant \frac{4}{(r + 1)^3}\chi
^2(P\vert \vert Q)
\end{equation}

\noindent and
\begin{align}
\frac{4}{(R + 1)^3}& \left[ {(R - 1)(1 - r) - \chi ^2(P\vert \vert
Q)} \right] \label{eq65}\\
& \leqslant \frac{2(R - 1)(1 - r)}{(R + 1)(1 + r)} - \Delta
(P\vert \vert Q)\notag\\
& \leqslant \frac{4}{(r + 1)^3}\left[ {(R-1)(1-r) - \chi ^2(P\vert
\vert Q)} \right].\notag
\end{align}
\end{proposition}

\begin{proof} Take $s = 2$ in (\ref{eq52}), (\ref{eq53})
and (\ref{eq54}) we get respectively (\ref{eq63}), (\ref{eq64})
and (\ref{eq65}).
\end{proof}

We observe that the Theorem \ref{the42} is not valid for $s = 0$,
$\frac{1}{2}$ and $1$. These particular values of $s$ we shall do
separately. In these cases, we don't have inequalities on both
sides as in the case of Propositions \ref{pro41} and \ref{pro42}.

\begin{proposition} \label{pro43} The following inequalities hold:
\begin{equation}
\label{eq66}
0 \leqslant \Delta (P\vert \vert Q) \leqslant \frac{32}{27}K(Q\vert \vert
P),
\end{equation}
\begin{equation}
\label{eq67}
0 \leqslant \Delta ^ * (P\vert \vert Q) \leqslant \frac{32}{27}\left[ {\chi
^2(Q\vert \vert P) - K(Q\vert \vert P)} \right]
\end{equation}

\noindent and
\begin{align}
0 & \leqslant \frac{32}{27}K(Q\vert \vert P) - \Delta (P\vert
\vert Q)\label{eq68}\\
&  \leqslant \frac{32}{27}\frac{(R - 1)\ln \frac{1}{r} + (1 -
r)\ln \frac{1}{R}}{R - r} - \frac{2(R - 1)(1 - r)}{(R + 1)(1 +
r)}\notag
\end{align}
\end{proposition}

\begin{proof} For $s = 0$ in (\ref{eq56}), we have
\begin{equation}
\label{eq69} g_\Delta (x) = \frac{8x^2}{(x + 1)^3}.
\end{equation}

This gives
\begin{equation}
\label{eq70}
{g}'_\Delta (x) = - \frac{8x(x - 2)}{(x + 1)^4}
\begin{cases}
 { \geqslant 0,} & {x \leqslant 2} \\
 { \leqslant 0,} & {x \geqslant 2} \\
\end{cases}.
\end{equation}

Thus we conclude that the function $g_W (x)$ given by (\ref{eq69})
is increasing in $x \in (0,2)$ and decreasing in $x \in (2,\infty
)$, and hence
\begin{equation}
\label{eq71}
M = \mathop {\sup }\limits_{x \in (0,\infty )} g_\Delta (x) = \mathop {\max
}\limits_{x \in (0,\infty )} g_\Delta (x) = g_\Delta (2) = \frac{32}{27}.
\end{equation}

Now (\ref{eq71}) together with (\ref{eq38}), (\ref{eq39}) and (\ref{eq40}) give respectively (\ref{eq66}),
(\ref{eq67}) and (\ref{eq68}).
\end{proof}

\begin{proposition} \label{pro44} The following inequalities hold:
\begin{equation}
\label{eq72}
0 \leqslant \Delta (P\vert \vert Q) \leqslant 4\mbox{ }h(P\vert \vert Q),
\end{equation}
\begin{equation}
\label{eq73}
0 \leqslant \Delta ^ * (P\vert \vert Q) \leqslant 2\mbox{ }\sum\limits_{i =
1}^n {(q_i - p_i )\sqrt {\frac{q_i }{p_i }} - 4\mbox{ }h(P\vert \vert Q)}
\end{equation}

\noindent and
\begin{align}
0 & \leqslant 4\mbox{ }h(P\vert \vert Q) - \Delta (P\vert \vert
Q)\label{eq74}\\
& \leqslant \frac{4(\sqrt R - 1)(1 - \sqrt r )}{\sqrt R + \sqrt r
} - \frac{2(R - 1)(1 - r)}{(R + 1)(1 + r)}.\notag
\end{align}
\end{proposition}

\begin{proof} For $s = \frac{1}{2}$ in (\ref{eq56}), we have
\begin{equation}
\label{eq75} g_\Delta (x) = \frac{8x^{3 / 2}}{(x + 1)^3}.
\end{equation}

This gives
\begin{equation}
\label{eq76}
{g}'_\Delta (x) = - \frac{12\sqrt x (x - 1)}{(x + 1)^4}
\begin{cases}
 { \geqslant 0,} & {x \leqslant 1} \\
 { \leqslant 0,} & {x \geqslant 1} \\
\end{cases}.
\end{equation}

Thus we conclude that the function $g_\Delta (x)$ given by
(\ref{eq75}) is increasing in $x \in (0,1)$ and decreasing in $x
\in (1,\infty )$, and hence
\begin{equation}
\label{eq77}
M = \mathop {\sup }\limits_{x \in (0,\infty )} g_\Delta (x) = \mathop {\max
}\limits_{x \in (0,\infty )} g_\Delta (x) = g_\Delta (1) = 1.
\end{equation}

Now (\ref{eq77}) together with (\ref{eq38}), (\ref{eq39}) and
(\ref{eq40}) give respectively (\ref{eq72}), (\ref{eq73}) and
(\ref{eq74}).
\end{proof}

\begin{proposition} \label{pro45} We have following inequalities:
\begin{equation}
\label{eq78} 0 \leqslant \Delta (P\vert \vert Q) \leqslant
\frac{32}{27}\, K(P\vert \vert Q),
\end{equation}
\begin{equation}
\label{eq79} 0 \leqslant \Delta ^ * (P\vert \vert Q) \leqslant
\frac{32}{27}\, K(Q\vert \vert P)
\end{equation}

\noindent and
\begin{align}
0 & \leqslant \frac{32}{27}\,K(P\vert \vert Q) - \Delta (P\vert
\vert Q)\label{eq80}\\
& \leqslant \frac{32}{27}\frac{(R - 1)r\ln r + (1 - r)R\ln R}{R -
r} - \frac{2(R - 1)(1 - r)}{(R + 1)(1 + r)}.\notag
\end{align}
\end{proposition}

\begin{proof} For $s = 1$ in (\ref{eq56}), we have
\begin{equation}
\label{eq81} g_\Delta (x) = \frac{8x}{(x + 1)^3}.
\end{equation}

This gives
\begin{equation}
\label{eq82}
{g}'_W (x) = - \frac{8(2x - 1)}{(x + 1)^4}
 = \begin{cases}
 { \geqslant 0,} & {x \leqslant \frac{1}{2}} \\
 { \leqslant 0,} & {x \geqslant \frac{1}{2}} \\
\end{cases}.
\end{equation}

Thus we conclude that the function $g_\Delta (x)$ given by
(\ref{eq81}) is increasing in $x \in (0,\frac{1}{2})$ and
decreasing in $x \in (\frac{1}{2},\infty )$, and hence
\begin{equation}
\label{eq83}
M = \mathop {\sup }\limits_{x \in (0,\infty )} g_\Delta (x) = \mathop {\max
}\limits_{x \in (0,\infty )} g_\Delta (x) = g(\frac{1}{2}) = \frac{32}{27}.
\end{equation}

Now (\ref{eq83}) together with (\ref{eq38}), (\ref{eq39}) and
(\ref{eq40}) give respectively (\ref{eq78}), (\ref{eq79}) and
(\ref{eq80}).
\end{proof}

\begin{remark} In view of relation (\ref{eq7}) and Propositions \ref{pro41}-\ref{pro45}, we have the following main bounds on
harmonic mean divergence:
\begin{equation}
\label{eq84} \frac{2r^3}{(r + 1)^3}\chi ^2(Q\vert \vert P)
\leqslant 1 - W(P\vert \vert Q) \leqslant \frac{2R^3}{(R +
1)^3}\chi ^2(Q\vert \vert P),
\end{equation}
\begin{equation}
\label{eq85} \frac{2}{(R + 1)^3}\chi ^2(P\vert \vert Q) \leqslant
1 - W(P\vert \vert Q) \leqslant \frac{2}{(r + 1)^3}\chi ^2(P\vert
\vert Q),
\end{equation}
\begin{equation}
\label{eq86} 0 \leqslant 1 - W(P\vert \vert Q) \leqslant
\frac{16}{27}K(Q\vert \vert P),
\end{equation}
\begin{equation}
\label{eq130} 0 \leqslant 1 - W(P\vert \vert Q) \leqslant 2\mbox{
}h(P\vert \vert Q)
\end{equation}

\noindent and
\begin{equation}
\label{131} 0 \leqslant 1 - W(P\vert \vert Q) \leqslant
\frac{16}{27}K(P\vert \vert Q).
\end{equation}
\end{remark}

The inequalities (\ref{eq85}) were also studied by Dragomir
\cite{dra3}. The inequalities (\ref{131}) can be seen in Dragomir
\cite{dra5}. The inequalities (\ref{eq72}) can be seen been in
Dragomir \cite{dra6} and Tops$\phi $e \cite{top}. The inequalities
(\ref{eq78}) can be in and Dragomir \cite{dra5}.

\section{Symmetric Chi-square Divergence and Inequalities}

In this section, we shall give bounds on \textit{symmetric
chi-square divergence} based on the Theorems \ref{the31} and
\ref{the32}.

\begin{theorem} \label{the51} For all $P,Q \in \Gamma _n $, we have
the following inequalities:
\begin{equation}
\label{eq87} 0 \leqslant \Psi(P\vert \vert Q) \leqslant \rho _\Psi
(P\vert \vert Q),
\end{equation}

\noindent where
\begin{equation}
\label{eq88} \rho _\Psi (P\vert \vert Q) = \Psi(P\vert \vert Q) +
\sum\limits_{i = 1}^n {\frac{(p_i - q_i )^2(p_i^2 + q_i^2 )}{p_i^2
q_i }} .
\end{equation}

If there exists $r,R$ $(0 < r \leqslant 1 \leqslant R < \infty )$
such that
\[
0 < r \leqslant \frac{p_i }{q_i } \leqslant R < \infty , \quad
\forall i \in \{1,2,...,n\},
\]

\noindent then we have the following inequalities:
\begin{equation}
\label{eq89} 0 \leqslant \Psi(P\vert \vert Q) \leqslant \alpha
_\Psi (r,R),
\end{equation}
\begin{equation}
\label{eq90} 0 \leqslant \Psi(P\vert \vert Q) \leqslant \beta
_\Psi (r,R)
\end{equation}

\noindent and
\begin{align}
0 & \leqslant \beta _\Psi (r,R) - \Psi(P\vert \vert Q)\label{eq91}\\
& \leqslant \gamma _\Psi (r,R)\left[ {(R - 1)(1 - r) - \chi
^2(P\vert \vert Q)} \right] \leqslant \alpha _\Psi (r,R),\notag
\end{align}

\noindent where
\begin{align}
\alpha _\Psi (r,R)&  = \frac{1}{4}(R - r)^2\left[ {2L_2^{ - 1}
(r,R) - L_1^{ - 1} (r,R)} \right],\label{eq92}\\
\beta _\Psi (r,R) & = (R - 1)(1 - r)(R + r)\label{eq93}\\
\intertext{and}
\label{eq94} \gamma _\Psi (r,R)& = 2L_2^{ - 1}
(r,R) - L_1^{ - 1} (r,R).
\end{align}
\end{theorem}

\begin{proof} Follows from Theorem \ref{the31} by
considering $f$ by $f_\Psi $ and making necessary calculations.
\end{proof}

\begin{theorem} \label{the52} Let $P,Q \in \Gamma _n $ and $s \in
\mathbb{R}$. Let there exists $r,R$ $(0 \leqslant r \leqslant 1
\leqslant R \leqslant \infty )$ such that $0 < r \leqslant
\frac{p_i }{q_i } \leqslant R < \infty $, $\forall i \in
\{1,2,...,n\}$.

$(a)$ For $s \leqslant - 1$, we have the following inequalities:
\begin{equation}
\label{eq95} \frac{2(r^3 + 1)}{r^{1 + s}}\Phi _s (P\vert \vert Q)
 \leqslant \Psi(P\vert \vert Q) \leqslant
\frac{2(R^3 + 1)}{R^{1 + s}}\Phi _s (P\vert \vert Q),
\end{equation}
\begin{align}
 \frac{2(r^3 + 1)}{r^{1 + s}}& \left[ {\rho _{\Phi _s }
(P\vert \vert Q) - \Phi _s (P\vert \vert Q)} \right]
\label{eq96}\\
& \leqslant \Psi^\ast (P\vert \vert Q)
 \leqslant \frac{2(R^3 + 1)}{R^{1 + s}}\left[ {\rho _{\Phi _s } (P\vert
\vert Q) - \Phi _s (P\vert \vert Q)} \right]\notag
\end{align}

\noindent and
\begin{align}
\frac{2(r^3 + 1)}{r^{1 + s}}& \left[ {\beta _{\Phi _s } (r,R) -
\Phi _s (P\vert \vert Q)} \right] \label{eq97}\\
& \leqslant \beta _\Psi (r,R) - \Psi(P\vert \vert Q)\notag\\
& \leqslant \frac{2(R^3 + 1)}{R^{1 + s}}\left[ {\beta _{\Phi _s }
(r,R) - \Phi _s (P\vert \vert Q)} \right].\notag
\end{align}

$(b)$ For $s \geqslant 2$, we have the following inequalities:
\begin{equation}
\label{eq98} \frac{2(R^3 + 1)}{R^{1 + s}}\Phi _s (P\vert \vert Q)
 \leqslant \Psi(P\vert \vert Q) \leqslant
\frac{2(r^3 + 1)}{r^{1 + s}}\Phi _s (P\vert \vert Q),
\end{equation}
\begin{align}
\frac{2(r^3 + 1)}{r^{1 + s}}& \left[ {\rho _{\Phi _s } (P\vert
\vert Q) - \Phi _s (P\vert \vert Q)} \right]
\label{eq99}\\
& \leqslant \Psi^\ast (P\vert \vert Q)
 \leqslant \frac{2(r^3 + 1)}{r^{1 + s}}\left[ {\rho _{\Phi _s } (P\vert
\vert Q) - \Phi _s (P\vert \vert Q)} \right] \notag
\end{align}

\noindent and
\begin{align}
\frac{2(R^3 + 1)}{R^{1 + s}}& \left[ {\beta _{\Phi _s } (r,R) -
\Phi _s (P\vert \vert Q)} \right]\label{eq100}\\
& \leqslant \beta _\Psi (r,R) - \Psi(P\vert \vert Q)\notag\\
& \leqslant \frac{2(r^3 + 1)}{r^{1 + s}}\left[ {\beta _{\Phi _s }
(r,R) - \Phi _s (P\vert \vert Q)} \right],\notag
\end{align}

\noindent where
\begin{equation}
\label{eq101} \Psi^\ast (P\vert \vert Q) = \rho _\Psi (P\vert
\vert Q) - \Psi(P\vert \vert Q) = \sum\limits_{i = 1}^n
{\frac{(p_i - q_i )^2(p_i^2 + q_i^2 )}{p_i^2 q_i }} .
\end{equation}
\end{theorem}

\begin{proof} Let us consider
\begin{equation}
\label{eq102} g_\Psi (x) = x^{2 - s}{f}''_\Psi (x) = 2x^{ - 1 -
s}(x^3 + 1), \quad x \in (0,\infty ),
\end{equation}

\noindent where ${f}''_\Psi (x)$ is as given by (\ref{eq18}).

We have
\begin{equation}
\label{eq103} {g}'_\Psi (x) = - 2x^{ - 2 - s}\left[ {(s - 2)x^3 +
(s + 1)} \right] \begin{cases}
 { \geqslant 0,} & {s \leqslant - 1} \\
 { \leqslant 0,} & {s \geqslant 2} \\
\end{cases}.
\end{equation}

From (\ref{eq103}), we conclude the followings:
\begin{equation}
\label{eq104} m = \mathop {\inf }\limits_{x \in [r,R]} g_\Psi (x)
= \mathop {\min }\limits_{x \in [r,R]} g_\Psi (x) = \begin{cases}
 {\frac{2(r^3 + 1)}{r^{1 + s}},} & {s \leqslant - 1} \\
 {\frac{2(R^3 + 1)}{R^{1 + s}},} & {s \geqslant 2} \\
\end{cases}
\end{equation}

\noindent and
\begin{equation}
\label{eq105} M = \mathop {\sup }\limits_{x \in [r,R]} g_\Psi (x)
= \mathop {\max }\limits_{x \in [r,R]} g_\Psi (x) = \begin{cases}
 {\frac{2(R^3 + 1)}{R^{1 + s}},} & {s \leqslant - 1} \\
 {\frac{2(r^3 + 1)}{r^{1 + s}},} & {s \geqslant 2} \\
\end{cases}
\end{equation}

In view of (\ref{eq104}) and (\ref{eq105}) and Theorem
\ref{the32}, we have the required proof.
\end{proof}

\begin{proposition} \label{51} We have the following bounds in
terms of $\chi ^2 - $divergence:
\begin{equation}
\label{eq106} (r^3 + 1)\chi ^2(Q\vert \vert P) \leqslant
\Psi(P\vert \vert Q) \leqslant (R^3 + 1)\chi ^2(Q\vert \vert P),
\end{equation}
\begin{align}
2(r^3 + 1)& \left[ {3\,\Phi _3 (Q\vert \vert P) - \chi
^2(Q\vert \vert P)} \right] \label{eq107}\\
& \leqslant \Psi^\ast (P\vert \vert Q)
 \leqslant 2(R^3 + 1)\left[ {3, \Phi _3 (P\vert \vert Q) - \chi
^2(Q\vert \vert P)} \right]\notag
\end{align}

\noindent and
\begin{align}
(r^3 + 1) & \left[ {\frac{(R - 1)(1 - r)}{rR} - \chi ^2(Q\vert
\vert P)} \right] \label{eq108}\\
& \leqslant (R - 1)(1 - r)(R + r) - \Psi(P\vert \vert Q)\notag\\
& \leqslant (R^3 + 1)\left[ {\frac{(R - 1)(1 - r)}{rR} - \chi
^2(Q\vert \vert P)} \right].\notag
\end{align}
\end{proposition}

\begin{proof} Take $s = - 1$ in (\ref{eq95}), (\ref{eq96})
and (\ref{eq97}) we get respectively (\ref{eq106}), (\ref{eq107})
and (\ref{eq108}).
\end{proof}

\begin{proposition} \label{pro52} The following bounds on
in terms of $\chi ^2 - $divergence hold:
\begin{equation}
\label{eq109} \frac{R^3 + 1}{R^3}\chi ^2(P\vert \vert Q) \leqslant
\Psi(P\vert \vert Q) \leqslant \frac{r^3 + 1}{r^3}\chi ^2(P\vert
\vert Q),
\end{equation}
\begin{equation}
\label{eq110} \frac{R^3 + 1}{R^3}\chi ^2(P\vert \vert Q) \leqslant
\Psi^\ast (P\vert \vert Q) \leqslant \frac{r^3 + 1}{r^3}\chi
^2(P\vert \vert Q)
\end{equation}

\noindent and
\begin{align}
\frac{R^3 + 1}{R^3}& \left[ {(R - 1)(1 - r) - \chi ^2(P\vert \vert
Q)} \right] \label{eq111}\\
& \leqslant \beta _\Psi (r,R) - \Psi(P\vert \vert Q)\notag\\
& \leqslant \frac{r^3 + 1}{r^3}\left[ {(R - 1)(1 - r) - \chi
^2(P\vert \vert Q)} \right].\notag
\end{align}
\end{proposition}

\begin{proof} Take $s = 2$ in (\ref{eq98}), (\ref{eq99}) and
(\ref{eq100}) we get respectively (\ref{eq109}), (\ref{eq110}) and
(\ref{eq111}).
\end{proof}

We observe that the above two propositions follows from Theorem
\ref{the52} immediately by taking $s = - 1$ and $s = 2$
respectively. But still there are another values of $s$ such as $s
= 0$, $s = 1$ and $s = \frac{1}{2}$ for which we can obtain
bounds. These values are studied below.

\begin{proposition} \label{pro53} We have following bounds in terms of
relative information:
\begin{equation}
\label{eq112} 0 \leqslant 3\sqrt[3]{2}\mbox{ }K(Q\vert \vert P)
\leqslant \Psi(P\vert \vert Q),
\end{equation}
\begin{equation}
\label{eq113} 0 \leqslant 3\sqrt[3]{2}\left[ {\chi ^2(Q\vert \vert
P) - K(Q\vert \vert P)} \right] \leqslant \Psi^\ast (P\vert \vert
Q)
\end{equation}

\noindent and
\begin{align}
0 & \leqslant \Psi(P\vert \vert Q) - 3\sqrt[3]{2}\, K(Q\vert \vert
P)\label{eq114}\\
& \leqslant (R - 1)(1 - r)(R + r) - 3\sqrt[3]{2}\, \frac{(R -
1)\ln \frac{1}{r} + (1 - r)\ln \frac{1}{R}}{R - r}.\notag
\end{align}
\end{proposition}

\begin{proof} For $s = 0$ in (\ref{eq102}), we have
\begin{equation}
\label{eq115} g_\Psi (x) = \frac{2(x^3 + 1)}{x},
\end{equation}

This gives
\begin{align}
{g}'_\Psi (x) & = \frac{2(2x^3 - 1)}{x^2}\label{eq116}\\
& = \frac{2(\sqrt[3]{2}\mbox{ }x - 1)(\sqrt[3]{4}\mbox{ }x^2 +
\sqrt[3]{2}\mbox{ }x + 1)}{x^2} \begin{cases}
 { \geqslant 0,} & {x \geqslant \frac{1}{\sqrt[3]{2}}} \\
 { \leqslant 0,} & {x \leqslant \frac{1}{\sqrt[3]{2}}} \\
\end{cases}.\notag
\end{align}

Thus we conclude that the function $g_\Psi (x)$ given by
(\ref{eq115}) is decreasing in $x \in (0,\frac{1}{\sqrt[3]{2}})$
and increasing in $x \in (\frac{1}{\sqrt[3]{2}},\infty)$, and
hence
\begin{equation}
\label{eq117} m = \mathop {\inf }\limits_{x \in (0,\infty )}
g_\Psi (x) = \mathop {\min }\limits_{x \in (0,\infty )} g_\Psi (x)
= g_\Psi (\frac{1}{\sqrt[3]{2}}) = 3\sqrt[3]{2}.
\end{equation}

Now (\ref{eq117}) together with (\ref{eq38}), (\ref{eq39}) and
(\ref{eq40}) give respectively (\ref{eq112}), (\ref{eq113}) and
(\ref{eq114}).
\end{proof}

\begin{proposition} \label{pro54} We have following bounds in terms of
Hellinger's discrimination:
\begin{equation}
\label{eq118} 0 \leqslant 16\mbox{ }h(P\vert \vert Q) \leqslant
\Psi(P\vert \vert Q),
\end{equation}
\begin{equation}
\label{eq119} 0 \leqslant 16\left[ {\frac{1}{2}\sum\limits_{i =
1}^n {(q_i - p_i )\sqrt {\frac{q_i }{p_i }} } - h(Q\vert \vert P)}
\right] \leqslant \Psi^\ast (P\vert \vert Q)
\end{equation}

\noindent and
\begin{align}
0 & \leqslant \Psi(P\vert \vert Q) - 16\mbox{ }h(P\vert \vert
Q)\label{eq120}\\
& \leqslant (R - 1)(1 - r)(R + r) - \frac{16(\sqrt R - 1)(1 -
\sqrt r )}{\sqrt R + \sqrt r }.\notag
\end{align}
\end{proposition}

\begin{proof} For $s = \frac{1}{2}$ in (\ref{eq102}), we
have
\begin{equation}
\label{eq121} g_\Psi (x) = \frac{2(x^3 + 1)}{x^{3 / 2}}.
\end{equation}

This gives
\begin{equation}
\label{eq122} {g}'_\Psi (x) = \frac{3(x^3 - 1)}{x^{5 / 2}} =
\frac{3(x - 1)(x^2 + x + 1)}{x^{5 / 2}} \begin{cases}
 { \geqslant 0,} & {x \geqslant 1} \\
 { \leqslant 0,} & {x \leqslant 1} \\
\end{cases}.
\end{equation}

Thus we conclude that the function $g_\Psi (x)$ given by
(\ref{eq121}) is decreasing in $x \in (0,1)$ and increasing in $x
\in (1,\infty )$, and hence
\begin{equation}
\label{eq123} m = \mathop {\inf }\limits_{x \in (0,\infty )}
g_\Psi (x) = \mathop {\min }\limits_{x \in (0,\infty )} g_\Psi (x)
= g_\Psi (1) = 4.
\end{equation}

Now (\ref{eq123}) together with (\ref{eq38}), (\ref{eq39}) and
(\ref{eq40}) give respectively (\ref{eq118}), (\ref{eq119}) and
(\ref{eq120}).
\end{proof}

\begin{proposition} \label{pro55} We have the following bounds in
terms of relative information:
\begin{equation}
\label{eq124} 0 \leqslant 3\sqrt[3]{2}\mbox{ }K(P\vert \vert Q)
\leqslant \Psi(P\vert \vert Q),
\end{equation}
\begin{equation}
\label{eq125} 0 \leqslant 3\sqrt[3]{2}\mbox{ }K(Q\vert \vert P)
\leqslant \Psi^\ast (P\vert \vert Q)
\end{equation}

\noindent and
\begin{align}
0 & \leqslant \Psi(P\vert \vert Q) - 3\sqrt[3]{2}\mbox{ }K(P\vert
\vert Q)\label{eq126}\\
& \leqslant (R - 1)(1 - r)(R + r) - 3\sqrt[3]{2}\mbox{ }\frac{(R -
1)r\ln r + (1 - r)R\ln R}{R - r}.\notag
\end{align}
\end{proposition}

\begin{proof} For $s = 1$ in (\ref{eq102}), we have
\begin{equation}
\label{eq127} g_\Psi (x) = \frac{2(x^3 + 1)}{x^2}.
\end{equation}

This gives
\begin{align}
{g}'_\Psi (x)& = \frac{2(x^3 - 2)}{x^3}\label{eq128}\\
&  = \frac{2(x - \sqrt[3]{2})(x^2 + \sqrt[3]{2}\mbox{ }x +
\sqrt[3]{4})}{x^3} \begin{cases}
 { \geqslant 0,} & {x \geqslant \sqrt[3]{2}\mbox{ }} \\
 { \leqslant 0,} & {x \leqslant \sqrt[3]{2}\mbox{ }} \\
\end{cases}.\notag
\end{align}

Thus we conclude that the function $g_\Psi (x)$ given by
(\ref{eq127}) is decreasing in $x \in (0,\sqrt[3]{2})$ and
increasing in $x \in (\sqrt[3]{2},\infty )$, and hence
\begin{equation}
\label{eq129} m = \mathop {\inf }\limits_{x \in (0,\infty )}
g_\Psi (x) = \mathop {\min }\limits_{x \in (0,\infty )} g_\Psi (x)
= g_\Psi (\sqrt[3]{2}) = 3\sqrt[3]{2}.
\end{equation}

Now (\ref{eq129}) together with (\ref{eq38}), (\ref{eq39}) and
(\ref{eq40}) give respectively (\ref{eq124}), (\ref{eq125}) and
(\ref{eq126}).
\end{proof}


\begin{thebibliography}{99}
\bibitem{pea} K. PEARSON, On the Criterion that a given system of deviations
from the probable in the case of correlated system of variables is
such that it can be reasonable supposed to have arisen from random
sampling, \textit{Phil. Mag.}, \textbf{50}(1900), 157-172.
\end{thebibliography}
\end{document}